\documentclass[a4paper,11pt]{amsart}
\usepackage{amssymb,amscd,verbatim, amsthm,amsmath,amsgen,xspace}
\usepackage{latexsym}
\usepackage{amsfonts}
\usepackage{color}
\usepackage[all]{xy}

\newcommand \ch[1]{{\check{#1}}}

\newcommand\Hom{\operatorname{Hom}}
\newcommand\Ad{\operatorname{Ad}}

\newcommand{\Kt}{\widetilde{K}}

\def\dim{\mathop{\hbox {dim}}\nolimits}
\def\Ad{\mathop{\hbox {Ad}}\nolimits}
\def\ad{\mathop{\hbox {ad}}\nolimits}

\def\ch{\mathop{\hbox {ch}}\nolimits}

\def\det{\mathop{\hbox {det}}\nolimits}

\def\Hom{\mathop{\hbox {Hom}}\nolimits}

\DeclareMathOperator\im{Im}

\def\ker{\mathop{\hbox{Ker}}\nolimits}
\DeclareMathOperator\Ker{Ker}

\def\R{\bbR}

\newcommand{\pf}{\begin{proof}}
\newcommand{\epf}{\end{proof}}
\newcommand{\eq}{\begin{equation}}
\newcommand{\eeq}{\end{equation}}
\newcommand{\eqn}{\begin{equation*}}
\newcommand{\eeqn}{\end{equation*}}

\newcommand{\frb}{\mathfrak{b}}

\newcommand{\frg}{\mathfrak{g}}
\newcommand{\frh}{\mathfrak{h}}
\newcommand{\frk}{\mathfrak{k}}
\newcommand{\frl}{\mathfrak{l}}

\newcommand{\frp}{\mathfrak{p}}
\newcommand{\frq}{\mathfrak{q}}
\newcommand{\frr}{\mathfrak{r}}
\newcommand{\frs}{\mathfrak{s}}
\newcommand{\frt}{\mathfrak{t}}

\newcommand{\frso}{\mathfrak{so}}

\renewcommand{\k}{{\mathfrak{k}}}
\renewcommand{\t}{{\mathfrak{t}}}

\newcommand{\bbC}{\mathbb{C}}

\newcommand{\bbR}{\mathbb{R}}

\newcommand{\caO}{\mathcal{O}}
\newcommand{\caS}{\mathcal{S}}

\newcommand{\tr}{\operatorname{tr}}

\newtheorem{thm}[equation]{Theorem}
\newtheorem{cor}[equation]{Corollary}
\newtheorem{lemma}[equation]{Lemma}
\newtheorem{prop}[equation]{Proposition}
\newtheorem{defi}[equation]{Definition}

\newtheorem{rmk}[equation]{Remark}

\numberwithin{equation}{section}

\let\ssize\scriptstyle
\newif\ifFIRST\newdimen\MAXright\MAXright0pt
\def\sdynkin{\bgroup\eightpoint\dynkin}
\def\endsdynkin{\enddynkin\egroup}
\def\dynkin{\bgroup\FIRSTtrue\hskip.5em\setbox1\hbox{$\diagup$}%
	\setbox2\hbox{$\diagdown$}%
	\setbox0\hbox to2\wd1{\hrulefill}%
	\setbox3\hbox{$\bullet$}%
	\setbox4\hbox{$\times$}%
	\setbox7\hbox{$\circ$}%       (L.K.)
	\def\whiteroot##1{\ifFIRST\setbox5\hbox{$##1$}\ifdim\wd5>1.3em%       (L.K.)
		\hskip-.5em\hskip.5\wd5\fi\fi\FIRSTfalse%                             (L.K.)
		\hskip-.25em\raise1.5\wd3\hbox to0pt{\hss\hskip.45em$%                (L.K.)
			\ssize##1$\hss}\copy7\hskip-.25em\setbox6\hbox{$##1$}%                (L.K.)
		\MAXright\wd6}%                                                       (L.K.)
	\def\root##1{\ifFIRST\setbox5\hbox{$##1$}\ifdim\wd5>1.3em%
		\hskip-.5em\hskip.5\wd5\fi\fi\FIRSTfalse%
		\hskip-.25em\raise1.5\wd3\hbox to0pt{\hss\hskip.45em$%
			\ssize##1$\hss}\copy3\hskip-.25em\setbox6\hbox{$##1$}%
		\MAXright\wd6}%
	\def\whitedroot##1{\ifFIRST\setbox5\hbox{$##1$}\ifdim\wd5>1.3em% (L.K.)
		\hskip-.5em\hskip.5\wd5\fi\fi\FIRSTfalse% (L.K.)
		\hskip-.25em\lower1.8\wd3\hbox to0pt{\hss\hskip.45em$%  (L.K.)
			\ssize##1$\hss}\copy7\hskip-.25em\setbox6\hbox{$##1$}% (L.K.)
		\MAXright\wd6}%
	\def\whiterroot##1{\hskip-.25em\copy7\hbox to0pt{\hskip.3em$\ssize##1$\hss}%
		\hskip-.25em\setbox6\hbox{\hskip.6em$##1##1$}%
		\MAXright\wd6}%
	\def\droot##1{\ifFIRST\setbox5\hbox{$##1$}\ifdim\wd5>1.3em%
		\hskip-.5em\hskip.5\wd5\fi\fi\FIRSTfalse%
		\hskip-.25em\lower1.8\wd3\hbox to0pt{\hss\hskip.45em$%
			\ssize##1$\hss}\copy3\hskip-.25em\setbox6\hbox{$##1$}%
		\MAXright\wd6}%
	\def\rroot##1{\hskip-.25em\copy3\hbox to0pt{\hskip.3em$\ssize##1$\hss}%
		\hskip-.25em\setbox6\hbox{\hskip.6em$##1##1$}%
		\MAXright\wd6}%
	\def\norroot##1{\hskip-.36em\copy4\hbox to0pt{\hskip.3em$\ssize##1$\hss}%
		\hskip-.48em\setbox6\hbox{\hskip.6em$##1##1$}%
		\MAXright\wd6}%
	\def\noroot##1{\ifFIRST\setbox5\hbox{$##1$}\ifdim\wd5>1.3em%
		\hskip-.5em\hskip.5\wd5\fi\fi\FIRSTfalse%
		\hskip-.36em\raise1.5\wd3\hbox to0pt{\hss\hskip.6em$%
			\ssize##1$\hss}\copy4\hskip-.38em\setbox6\hbox{$##1$}%
		\MAXright\wd6}%
	\def\nodroot##1{\ifFIRST\setbox5\hbox{$##1$}\ifdim\wd5>1.3em%
		\hskip-.5em\hskip.5\wd5\fi\fi\FIRSTfalse%
		\hskip-.36em\lower1.8\wd3\hbox to0pt{\hss\hskip.6em$%
			\ssize##1$\hss}\copy4\hskip-.38em\setbox6\hbox{$##1$}%
		\MAXright\wd6}%
	\def\nolink{\hskip\wd0}%      (L.K.)
	\def\link{\raise.22em\copy0}%
	\def\llink##1{\raise.32em\copy0\hskip-\wd0%
		\raise.12em\copy0\hskip-.5\wd0\hbox to0pt{\hss$##1$\hss}\hskip.5\wd0}%
	\def\lllink##1{\raise.22em\copy0\hskip-\wd0\raise.32em\copy0\hskip-\wd0%
		\raise.12em\copy0\hskip-.5\wd0\hbox to0pt{\hss$##1$\hss}\hskip.5\wd0}%
	\def\rootupright##1{\hbox to0pt{\raise.45em\copy1\hskip-.25em\raise1.3\ht1%
			\hbox{\copy3\hskip.3em$\ssize##1$}\hss}%
		\setbox6\hbox{\hskip.6em\copy1\copy1$##1##1$}%
		\ifdim\MAXright<\wd6\MAXright\wd6\fi}%
	\def\whiterootupright##1{\hbox to0pt{\raise.45em\copy1\hskip-.25em\raise1.3\ht1% (L.K.)
			\hbox{\copy7\hskip.3em$\ssize##1$}\hss}% (L.K.)
		\setbox6\hbox{\hskip.6em\copy1\copy1$##1##1$}% (L.K.)
		\ifdim\MAXright<\wd6\MAXright\wd6\fi}% (L.K.)
	\def\norootupright##1{\hbox to0pt{\raise.45em\copy1\hskip-.36em\raise1.3\ht1%
			\hbox{\copy4\hskip.3em$\ssize##1$}\hss}%
		\setbox6\hbox{\hskip.6em\copy1\copy1$##1##1$}%
		\ifdim\MAXright<\wd6\MAXright\wd6\fi}%
	\def\rootdownright##1{\hbox to0pt{\raise-.5em\copy2\hskip-.25em\raise-1.35\ht1%
			\hbox{\copy3\hskip.3em$\ssize##1$}\hss}\setbox6%
		\hbox{\hskip.6em\copy2\copy2$##1##1$}%
		\ifdim\MAXright<\wd6\MAXright\wd6\fi}%
	\def\whiterootdownright##1{\hbox to0pt{\raise-.5em\copy2\hskip-.25em\raise-1.35\ht1% (L.K.)
			\hbox{\copy7\hskip.3em$\ssize##1$}\hss}\setbox6% (L.K.)
		\hbox{\hskip.6em\copy2\copy2$##1##1$}% (L.K.)
		\ifdim\MAXright<\wd6\MAXright\wd6\fi}% (L.K.)
	\def\rootdown##1{\hbox to0pt{\hskip-.05em\vrule height.25em depth.65em%
			\hskip-.25em\raise-.95em\hbox{\copy3\hskip.3em$\ssize##1$}\hss}%
		\setbox6\hbox{$##1$}%
		\ifdim\MAXright<\wd6\MAXright\wd6\fi}%
	\def\whiterootdown##1{\hbox to0pt{\hskip-.05em\vrule height.25em depth.65em% (L.K.)
			\hskip-.25em\raise-.95em\hbox{\copy7\hskip.3em$\ssize##1$}\hss}% (L.K.)
		\setbox6\hbox{$##1$}% (L.K.)
		\ifdim\MAXright<\wd6\MAXright\wd6\fi}% (L.K.)
	\def\dots{\hskip.5em\cdots\hskip.5em}}%
\def\enddynkin{\ifdim\MAXright>1em\hskip.5\MAXright\else\hskip.5em\fi\egroup}%

\begin{document}
%\today

\title{Dirac cohomology and character lifting}

\author{Jing-Song Huang}

\address{Department of Mathematics, Hong Kong University of Science and Technology,
Clear Water Bay, Kowloon, Hong Kong SAR, China}
\email{mahuang@ust.hk}

\thanks{The research described in this paper is supported by grants No. 16303218  from
Research Grant Council of HKSAR}
	
\begin{abstract}
The endoscopic transfer factor is expressed as  
difference of characters for the even and odd parts of the spin modules, or Dirac index of the  trivial representation. The
lifting of tempered characters  in terms of
 index of Dirac cohomology is calculated explicitly.
\end{abstract}

%\dedicatory{Dedicated to the memory of Kostant}
\maketitle

%%%%%%%%%%%%%%  Section 1 %%%%%%%%%%%%%%
\section{Introduction}\label{section intro}

The lifting of representations between reductive algebraic groups plays an important role in representation theory.
The endoscopic transfer in Langlands functionality and the theta correspondence in Howe's reductive dual pairs are primary examples. 

Dirac operators are employed for geometric construction of discrete series by Parthasarathy \cite{P}, Atiyah and Schmid \cite{AS}, and tempered representations by Wolf \cite{W}. In the late 1990's, Vogan made a conjecture on the algebraic property of the Dirac operators in Lie algebra setting.  This conjecture was proved by Pandzic and myself in 2002 \cite{HP1}.  This led us to study Dirac cohomology of Harish-Chandra modules \cite{HP2}. 
Kostant  extended the concept of Dirac cohomology and Vogan's conjecture to the more general setting of the cubic Dirac operator \cite{Ko2}. 

In the formulation of central problems in Langlands program stable conjugacy plays an pivotal role. The theory of endoscopy investigates the difference between orbital integral over ordinary and stable conjugacy classes.  As Dirac cohomology of a Harish-Chandra module determines its K-character \cite{HPZ}, it corresponds to the dual object of the orbital integral on elliptic elements.  In this note we focus on using Dirac cohomology 
for endoscopic transfer.  The transfer factor
is difference of characters for the even and odd parts of the spin modules, or Dirac index of the  trivial representation \cite{H}. 
The aim of this note is to extend the calculation of lifting of characters to tempered representations with nonzero Dirac cohomology.
The ultimate goal is to understand lifting of characters for all unitary representations with nonzero Dirac cohomology (the Dirac series).

Jeff Adams \cite{A} defined a lifting on characters between orthogonal groups and nonlinear metaplectic groups over real numbers,
 and it was extended to the p-adic case by Tatiana Howard \cite{Hd}.
  This lifting of characters is closely related to both endoscopy and theta correspondence, 
which also appears in the work of David Renard \cite{R}
 and Wen-Wei Li \cite {Li}.
The Adams lifting of characters for orthogonal and
symplectic groups is given by the formal difference of the oscillator representations of 
metaplectic groups and are related to the symplectic Dirac cohomology for Lie superalgebra \cite{HP3}.
Since the symplectic Dirac cohomology of modules for Lie superalgebras is very different from
the Dirac cohomology of Harish-Chandra modules, we deal with the transfer factor of the Adams lifting
in another paper \cite{H2}.

%%%%%%%%%%%%%% Section 2  %%%%%%%%%%
\section{Preliminaries on Dirac cohomology}
%%%%%%%%%%%%%%%%%%%%%%%%%%%%%%%%%%%

For a real reductive group $G$ with a Cartan involution $\theta$, denote by $\frg_0$ its Lie algebra
and assume that $K=G^\theta$ is a maximal compact subgroup of $G$.
Let $\frg=\frk\oplus\frp$ be the Cartan decomposition
for the complexified Lie algebra of $G$.  Let $B$ be a non-degenerate invariant symmetric bilinear
form on $\frg$, which restricts to the Killing form on the semisimple part $[\frg,\frg]$ of $\frg$.

Let $U(\frg)$ be the universal enveloping
algebra of $\frg$ and $C(\frp)$ the Clifford algebra of
$\frp$ with respect to $B$. Then one can consider the following version of the Dirac
operator:
$$
D=\sum_{i=1}^n Z_i\otimes Z_i \in U(\frg)\otimes C(\frp);
$$
here $Z_1,\dots,Z_n$ is an orthonormal basis of $\frp$ with respect to the symmetric bilinear
form $B$.  It follows that $D$ is independent of the choice of the orthonomal basis $Z_1,\dots,Z_n$
and it is invariant under the diagonal adjoint action of $K$.

The Dirac operator $D$ is a square root of Laplace operator
associated to the symmetric pair $(\frg,\frk)$. To explain this, we
start with a Lie algebra map
\begin{equation*}
\alpha:\frk\rightarrow C(\frp)
\end{equation*}
which is defined by the adjoint map $\ad:\frk\rightarrow\frso(\frp)$
composed with the embedding of $\frso(\frp)$ into $C(\frp)$ using
the identification $\frso(\frp)\simeq\bigwedge^2\frp$. The explicit
formula for $\alpha$ is (see \cite[\S2.3.3]{HP2})
\begin{equation}
\alpha(X)=-\frac{1}{4}\sum_{j}[X,Z_j]Z_j.
\end{equation}
Using $\alpha$ we can embed the Lie algebra $\frk$ diagonally into
$U(\frg)\otimes C(\frp)$, by
\begin{equation*}\label{Delta map}
X\mapsto X_\Delta=X\otimes1+1\otimes\alpha(X).
\end{equation*}
This embedding extends to $U(\frk)$. We denote the image of $\frk$
by $\frk_\Delta$, and then the image of $U(\frk)$ is the enveloping
algebra $U(\frk_\Delta)$ of $\frk_\Delta$.

Let $\Omega_\frak g$ be the Casimir operator
for $\frak g$, given by $\Omega_\frak g = \sum Z_i^2 -\sum W_j^2$,
where $W_j$ is an orthonormal basis for $\frak k_0$ with respect to
the inner product $-B$, where $B$ is the Killing form. Let
$\Omega_{\frk}=-\sum W_j^2$ be the Casimir operator for $\frak k$.
The image of $\Omega_\frk$ under $\Delta$ is denoted by
$\Omega_{\frk_\Delta}$.

Then
\begin{equation}\label{D^2}
D^2 = -\Omega_\frak g\otimes 1 + \Omega_{\frak k_\Delta} +
(||\rho_c||^2-||\rho||^2)1\otimes 1,
\end{equation}
where $\rho$ and $\rho_c$ are half sums of positive roots and compact
positive roots respectively.

The Vogan conjecture says that every element $z\otimes 1$ of
$Z(\frg)\otimes 1 \subset U(\frg)\otimes C(\frp)$ can be written as
$$\zeta(z)+Da+bD$$
where $\zeta(z)$ is in $Z(\k_\Delta)$,
and $a,b\in U(\frg)\otimes C(\frp)$.

A main result in \cite{HP1}  is  introducing a differential $d$
on the $K$-invariants in $U(\frg)\otimes C(\frp)$ defined by a super bracket with $D$,
and determination of the cohomology of this
differential complex.  As a consequence,
Pand\v{z}i\'c and I proved the following theorem.
In the following we denote by $\frh$ a Cartan subalgebra of $\frg$
containing a Cartan subalgebra $\frt$ of $\frk$ so that $\frt^*$ is embedded into $\frh^*$,
and by $W$ and $W_K$ the Weyl groups of $(\frg,\frh)$ and $(\frk,\frt)$
respectively.

\begin{thm}[\cite{HP1}]
Let $\zeta:
Z(\frg)\rightarrow Z(\frk)\cong Z(\frk_{\Delta})$ be the algebra homomorphism
that is determined
by the following commutative diagram:
\begin{equation*}\label{Vogan's diagram}
\CD
  Z(\mathfrak{g}) @> \zeta >> Z(\mathfrak{k}) \\
  @V \eta VV @V \eta_{\mathfrak{k}} VV  \\
  P(\mathfrak{h}^*)^{W} @>\mathrm{Res}>>
  P(\mathfrak{t}^*)^{W_K},
\endCD
\end{equation*}
where $P$ denotes the polynomial algebra, and  vertical maps $\eta$ and $\eta_\frk$  are Harish-Chandra
isomorphisms.
Then for each $z\in Z(\frg)$ one has
\begin{equation*}\label{Vogan's equation}
z\otimes 1-\zeta(z)=Da+aD, \text{\ for some\ }a\in U(\frg)\otimes C(\frp).
\end{equation*}

\end{thm}

For any admissible $(\frg,K)$-module $X$,
Vogan (\cite{V}, \cite{HP1}) introduced the notion of \emph{Dirac cohomology} $H_D(X)$ of $X$.
Consider  the action of the Dirac operator $D$  on $X\otimes S$,
with $S$ the spinor module  for the Clifford algebra $C(\frp)$. The Dirac cohomology is defined as follows:
$$H_D(X)\colon =\Ker D/ (\im D \cap \Ker D).$$
It follows from the identity (\ref{D^2}) that $H_D(X)$ is a finite-dimensional module
for the spin double cover $\Kt$ of $K$. In case $X$ is unitary, $H_D(X)=\ker D=\ker D^2$
since $D$ is self-adjoint with respect to a natural Hermitian inner product on $X\otimes S$.
As a consequence of the above theorem, we have that
 $H_D(X)$, if nonzero, determines the infinitesimal character of $X$.

\begin{thm}[\cite{HP1}]
Let $X$ be an admissible $(\frg,K)$-module with standard infinitesimal character parameter
$\Lambda\in \frh^*$.  Suppose that $H_D(X)$ contains a representation of
$\Kt$ with infinitesimal character $\lambda$.  Then
$\Lambda$ and $\lambda\in\frt^*\subseteq \frh^*$ are conjugate under $W$.
\end{thm}

The above theorem is proved in \cite{HP1} for a connected semisimple
Lie group $G$.  It is straightforward to extend the result to a
possibly disconnected reductive Lie group  in Harish-Chandra's class \cite{DH}.

Vogan's conjecture implies a refinement of the celebrated
Parthasarathy's Dirac inequality, which is an extremely useful tool
for the classification of irreducible unitary representations of
reductive Lie groups.

\begin{thm}[Extended Dirac Inequality \cite{P}, \cite{HP1}]
%\cite{VZ})}
Let $X$ be an irreducible unitary $(\frg,K)$-module with infinitesimal
character $\Lambda$.  Fix a representation of $K$ occurring in $X$
with a highest weight $\mu\in\t^*$, and a positive root system $\Delta^+(\frg)$
for $\frt$ in $\frg$.  Here $\t$ is a Cartan subalgebra of $\k$.
Write
$$\rho_c=\rho(\Delta^+(\frk)),\ \rho_n=\rho(\Delta^+(\frp)).$$
Fix an element $w\in W_K$ such that $w(\mu-\rho_n)$ is dominant
for $\Delta^+(\frk)$. Then
$$\langle w(\mu-\rho_n)+\rho_c,w(\mu-\rho_n)+\rho_c\rangle
\geq \langle \Lambda,\Lambda\rangle.$$
The equality  holds if and only if there exists a $w\in W$
such that
$$\Lambda = w(\mu-\rho_n)+\rho_c.$$
\end{thm}

%%%%%%%%%%%% Section 3  %%%%%%%%%%%%%%
\section{Dirac series and elliptic representations}
%%%%%%%%%%%%%%%%%%%%%%%%%%%%%%%%%%%

We say that an irreducible representation is in the {\it Dirac series} if it is a unitary representation with nonzero Dirac cohomology.
The Dirac series contains several families of very important unitary
representations including discrete series, unitary highest weight modules and unitary representations with nonzero
$(\frg,K)$-cohomology.  As shown in the previous section, a unitary $A_\frq(\lambda)$-module 
with admissible $\lambda$ is a Dirac series.  Any irreducible unitary highest
module is also a Dirac series. Thus, there exists Dirac series other than $A_\frq(\lambda)$-module.  
We now discuss the relationship between Dirac series and elliptic representations.

Let $G$ be a connected semisimple algebraic group
over a local field $F$ of characteristic $0$.  Arthur \cite{A1} studied a subset $\Pi_{temp, ell}(G(F))$
of tempered representations of $G(F)$, namely elliptic tempered representations.
The set of tempered representations $\Pi_{temp}(G(F))$ includes the discrete series
and in general the irreducible constituents of representations induced from discrete series.
These are exactly the representations which occur
in the Plancherel formula for $G(F)$.
 In Harish-Chandra's theory\cite{HC1,HC2}, the character of an infinite dimensional representation $\pi$ is defined as a distribution
$$\Theta(\pi,f)=\tr \big{(} \int_{G(F)}f(x)\pi(x)dx\ \big{)}, \  \ \ \  f\in C^\infty_c(G(F) ),$$
which can be identified with a function on $G(F)$.  In other words,
$$\Theta(\pi,f)=\int_{G(F)}f(x)\Theta(\pi,x)dx, \ \ \ \ f\in C^\infty_c(G(F)),$$
where $\Theta(\pi,x)$ is a locally integrable function on $G(F)$ that is smooth on the
open dense subset $G_{reg}(F)$ of regular elements.
A representation $\pi$ is called elliptic if $\Theta(\pi,x)$ does not vanish on the set of
elliptic elements in $G_{reg}(F)$.

The central objects in \cite{A1} are the normalized
characters $\Phi(\pi,\gamma)$, namely the functions defined by
$$\Phi(\pi,\gamma)=|D(\gamma)|^{1\over 2}\Theta(\pi,\gamma),\ \pi  \in \Pi_{temp, ell}(G(F)), \
\gamma\in G_{reg}(F),$$
where
$$D(\gamma)=\det(1-\Ad(\gamma))_{\frg/\frg_\gamma},$$
is the Weyl discriminant.  We will show how this normalized
character $\Phi(\pi,\gamma)$ is related to the Dirac cohomology of the
Harish-Chandra module of $\pi$
for a real group $G(\R)$.

From now on we are concerned only with the real group $G(\R)$.
Note that $G(\R)$ has elliptic elements if and only if
it is of equal rank with $K(\R)$.  We also assume
this equal rank condition.
Induced representations from proper parabolic subgroups are not elliptic.
Consider the quotient of Grothendieck group of the category of finite length
Harish-Chandra modules by the subspace generated by induced
representations.  Let us call this quotient group the elliptic Grothendieck
group.  Authur \cite{A1} found an orthonormal basis of this elliptic
Grothendieck group in terms of elliptic tempered (possibly virtual) characters.
Those characters are the super tempered distributions defined by Harish-Chandra \cite{HC3}.

The tempered elliptic representations for
the real group $G(\R)$ are the representations
with non-zero Dirac index, which are studied in \cite{Lab1}.
Labesse shows that the tempered elliptic representations are precisely the fundamental series.
We now discuss the general elliptic representations
and and their Dirac index.

Recall that if $X$ is an
admissible $(\frg,K)$-module with $K$-type decomposition
$X=\bigoplus_{\lambda}m_\lambda E_\lambda$, then the $K$-character of $X$ is the formal series
\[
\ch X=\sum_{\lambda}m_\lambda\ch E_\lambda,
\]
where $\ch E_\lambda$ is the character of the irreducible $K$-module $E_\lambda$.
Moreover, this definition makes sense also for virtual $(\frg,K)$-modules $X$;
in that case, the integers $m_\lambda$ can be negative. In the following
we will often deal with representations of the spin double cover $\Kt$ of $K$,
and not $K$, but we will still denote the corresponding character by $\ch$.

Since $\frp$ is even-dimensional, the spin module $S$ decomposes as $S^+\oplus S^-$, with the $\frk$-submodules
$S^\pm$ being the even respectively odd part of $S\cong{\textstyle\bigwedge}\frp^+$.
Let $X=X_\pi$ be the Harish-Chandra module of an irreducible admissible
representation $\pi$ of $G(\bbR)$.  We consider the following difference
of $\Kt$-modules, the spinor index of $X$:
\[
I(X)=X\otimes S^+-X\otimes S^-.
\]
It is a virtual $\Kt$-module, an integer combination of
finitely many $\Kt$-modules.
The Dirac operator $D$ induces the action of the following $\Kt$-equivariant operators
\[
D^{\pm}:X\otimes S^{\pm}\rightarrow X\otimes S^{\mp}.
\]
Since $D^2$ acts by a scalar on each $\Kt$-type, most of $\Kt$-modules in
$X\otimes S^+$ are the same as in $X\otimes S^-$.  It is straightforward to show the following identity.

\begin{lemma}\cite[Lemma 8.1]{H} The spinor index is equal to the Euler characteristic
of Dirac cohomology, i.e.,
\[
I(X)=H_D^+(X)-H_D^-(X).
\]
\end{lemma}

The spinor index $I(X)$ is also called the Dirac index of $X$,
since it is equal to the index of $D^+$, in the sense of index for
a Fredholm operator.  It is also identical to the Euler characteristic of Dirac cohomology $H_D(X)$.
We denote by $\theta(X)$ the character of $I(X)$.
In terms of characters, this reads
$$
\theta(X)=\ch I(X)=\ch X(\ch S^+-\ch S^-)=\ch H_D^+(X)-\ch H_D^-(X).
$$

If we view $\ch E_\lambda$ as functions on $K$, then  the series
$$\ch X=\sum_{\lambda}m_\lambda\ch E_\lambda$$
converges to a distribution on $K$ and it coincides with
$\Theta(X)$ on $K\cap  G_{reg}$, according to Harish-Chandra \cite{HC1}.
Then the absolute value $|\theta_\pi|$ coincides with the absolute value
$
|\Phi(\pi,\gamma)|=|D(\gamma)|^{1\over 2}|\Theta(\pi,\gamma)|
$
on regular elliptic elements.  We write this fact as the following lemma.

\begin{lemma} \label{character-theta}  For any regular elliptic elements $\gamma$, we have
$$|\theta_\pi(\gamma)|=|\Phi(\pi,\gamma)|.$$
\end{lemma}

The following theorem follows immediately from the lemma.

\begin{thm}\cite[Theorem 8.3]{H}
Let $\pi$ be an irreducible admissible representation of $G(\R)$
with Harish-Chandra module $X_\pi$.
Then $\pi$ is elliptic if and only if the Dirac index $I(X_\pi)\neq 0$.
\end{thm}

We also recall a result from \cite{H}.
\begin{thm}\cite[Theorem 10.5]{H}  Suppose $\pi$ is an irreducible unitary elliptic representation of
$G(\R)$ with a regular infinitesimal character.
Then $X_\pi\cong A_\frq(\lambda)$.
\end{thm}

As a consequence of the above two theorems, we have the following
\begin{cor} Suppose that rank of $G(\R)$ is equal to rank of $K(\R)$.  Then any Dirac series of $G(\R)$ with regular infinitesimal character is
an $A_\frq(\lambda)$-module.
\end{cor}

The Dirac index of a representation determines its character on compact Cartan subgroups.
As shown by Harish-Chandra \cite{HC3}, the character of a discrete series is determined 
completely on the set of regular elliptic elements. 
It is a natural question whether the $A_\frq(\lambda)$ has the same property,
namely whether the Dirac index determines the representation.
It was shown in \cite{HPV} that it is indeed true for most of simple Lie groups
except for a few exception.  In those exceptions, one needs the Dirac cohomology
together with the rank of $[\frl,\frl]$ to determine the corresponding representation $A_\frq(\lambda)$.

%%%%%%%%%%%    Section 4    %%%%%%%%%%%%%
\section{Cubic Dirac operators and associated cohomology}
%%%%%%%%%%%%%%%%%%%%%%%%%%%%%%%%%%%

We now recall the definition of Kostant's cubic Dirac operator
and the basic properties of the corresponding Dirac cohomology.
Let $\frg$ be a semisimple complex Lie algebra with Killing form $B$.
Let $\frr\subset\frg$ be a reductive Lie subalgebra such that
$B|_{\frr\times \frr}$  is non-degenerate. Let
$\frg=\frr\oplus\frs$ be the orthogonal decomposition with respect
to $B$. Then the restriction
$B|_\frs$ is also non-degenerate. Denote by $C(\frs)$ the Clifford
algebra of $\frs$ with
\begin{equation*}
uu'+u'u=-2B(u, u')
\end{equation*}
for all $u, u'\in\frs$. The above choice of sign is the same as in \cite{HP2}, but
different from the definition in \cite{Ko1}, as well as in \cite{Ko2}.  The two different
choices of signs have no essential difference
since  the two bilinear forms are equivalent over $\mathbb{C}$. Now
fix an orthonormal basis $Z_1, \ldots, Z_m$ of $\frs$. Kostant
\cite{Ko1} defines the cubic Dirac operator $D$ by
\begin{equation*}
D=\sum_{i=1}^m{Z_i\otimes Z_i+1\otimes v}\in U(\frg)\otimes C(\frs).
\end{equation*}
Here $v\in C(\frs)$ is the image of the fundamental 3-form
$w\in\bigwedge^3(\frs^*)$,
\begin{equation*}
w(X,Y,Z)=\frac{1}{2}B(X,[Y,Z]),
\end{equation*}
under the Chevalley map $\bigwedge(\frs^*)\rightarrow C(\frs)$ and
the identification of $\frs^*$ with $\frs$ by the Killing form $B$.
Explicitly,
\begin{equation*}
v=\frac{1}{2}\sum_{1\leq i<j<k \leq m}B([Z_i, Z_j],Z_k)Z_iZ_jZ_k.
\end{equation*}

The cubic Dirac operator has a good square in analogue with
the Dirac operator associated with the symmetric pair $(\frg,\frk)$ in Section 2.  We
have  a similar Lie algebra map
\begin{equation*}
\alpha:\frr\rightarrow C(\frs)
\end{equation*}
which is defined by the adjoint map $\ad:\frr\rightarrow\frso(\frs)$
composed with the embedding of $\frso(\frs)$ into $C(\frs)$ using
the identification $\frso(\frs)\simeq\bigwedge^2\frs$. The explicit
formula for $\alpha$ is (see \cite[\S 2.3.3]{HP2})
\begin{equation}\label{mapalpha}\alpha(X)=-\frac{1}{4}\sum_{j}[X, Z_j]Z_j, \quad\
X\in\frr.
\end{equation}
Using $\alpha$ we can embed the Lie algebra $\frr$ diagonally into
$U(\frg)\otimes C(\frs)$, by
\begin{equation*}\label{Delta map}
X\mapsto X_\Delta=X\otimes1+1\otimes\alpha(X).
\end{equation*}
This embedding extends to $U(\frr)$. We denote the image of $\frr$
by $\frr_\Delta$, and then the image of $U(\frr)$ is the enveloping
algebra $U(\frr_\Delta)$ of $\frr_\Delta$. Let $\Omega_\frg$ (resp.
$\Omega_\frr$) be the Casimir elements for $\frg$ (resp. $\frr$).
The image of $\Omega_\frr$ under $\Delta$ is denoted by
$\Omega_{\frr_\Delta}$.

Let $\frh_\frr$ be a Cartan subalgebra of $\frr$ which is contained
in $\frh$. It follows from Kostant's calculation (\cite{Ko1}, Theorem 2.16) that
\begin{equation}\label{square}
D^2=-\Omega_\frg\otimes1+\Omega_{\frr_\Delta}-(\|\rho\|^2 - \|\rho_{\frr}\|^2)1\otimes1,
\end{equation}
where $\rho_\frr$ denote the half sum of positive roots for $(\frr,
\frh_\frr)$.  We also note the sign difference with Kostant's formula
due to our choice of bilinear form for the definition of the Clifford algebra $C(\frs)$.

We denote by $W$ the Weyl group associated to the root system
$\Delta(\frg,\frh)$ and $W_\frr$ the Weyl group associated to the root system
$\Delta(\frr,\frh_\frr)$. The following theorem due to Kostant is an
extension of Vogan's conjecture on the symmetric pair case which
is proved in \cite{HP1}.
(See \cite{Ko2} Theorems 4.1 and 4.2 or \cite{HP2} Theorem 4.1.4).

\begin{thm}\label{Vogan's conjecture1}
There is an algebra homomorphism 
$$\zeta:
Z(\frg)\rightarrow Z(\frr)\cong Z(\frr_{\Delta})$$ such that for any $z\in Z(\frg)$ one has
\begin{equation*}\label{Vogan's equation}
z\otimes 1-\zeta(z)=Da+aD \text{\ for some\ }a\in U(\frg)\otimes C(\frs).
\end{equation*}
Moreover, $\zeta$ is determined
by the following commutative diagram:
\begin{equation*}\label{Vogan's diagram}
\CD
  Z(\mathfrak{g}) @> \zeta >> Z(\mathfrak{r}) \\
  @V \eta VV @V \eta_{\mathfrak{r}} VV  \\
  P(\mathfrak{h}^*)^{W} @>\mathrm{Res}>>
  P(\mathfrak{h}_\frr^*)^{W_{\mathfrak{r}}}.
\endCD
\end{equation*}
Here the vertical maps $\eta$ and $\eta_\frr$ are Harish-Chandra
isomorphisms.
\end{thm}

\begin{defi} Let $S$ be a spin module of $C(\frs)$.
Consider the action of $D$ on $V\otimes S$
\begin{equation}\label{Dirac map}
D:V\otimes S\rightarrow V\otimes S
\end{equation}
with $\frg$ acting on $V$ and $C(\frs)$ on $S$.  The {\it Dirac
cohomology} of $V$ is defined to be the $\frr$-module
\begin{equation*}
H_D(V):=\Ker D/(\Ker D\cap \im D).
\end{equation*}
\end{defi}

The following theorem is a consequence of the above theorem.
\begin{thm}[\cite{Ko2},\cite{HP2}]
Let $V$ be a $\frg$-module with
$Z(\frg)$ infinitesimal character $\chi_\Lambda$. Suppose that an
$\frr$-module $N$ is contained in the Dirac cohomology
$H_D(V)$ and has $Z(\frr)$ infinitesimal
character $\chi_\lambda$ . Then $\lambda=w\Lambda$ for some $w\in W$.
\end{thm}

Suppose that $V_\lambda$ is a finite-dimensional representation
with highest weight $\lambda\in \frh^*$.  Kostant \cite{Ko2} calculated
the Dirac cohomology of $V_\lambda$ with respect to any equal rank
quadratic subalgebra $\frr$ of $\frg$.  Assume that $\frh\subset\frr\subset\frg$ is
the Cartan subalgebra for both $\frr$ and $\frg$.
Define $W(\frg,\frh)^1$ to be the subset of the Weyl group $W(\frg,\frh)$ by
$$
W(\frg,\frh)^1=\{w\in W(\frg,\frh)\ | w(\rho) \text{ is }\Delta^+(\frr,\frh)-\text{dominant}
\}.
$$
This is the same as the subset of elements $w\in W(\frg,\frh)$ that map the positive Weyl $\frg$-chamber
into the positive $\frr$-chamber.  There is a bijection
$W(\frr,\frh)\times W(\frg,\frh)^1\rightarrow W(\frg,\frh)$ given by $(w,\tau)\mapsto w\tau$.
Kostant \cite{Ko2} proved the following result.

\begin{prop}[Kostant \cite{Ko2}]  One has 
$$
H_D(V_\lambda)=\bigoplus_{w\in W(\frg,\frh)^1}E_{w(\lambda+\rho)-\rho_\frr}.
$$
\end{prop}
%%%%%%%%%%%%%%   Section 5 %%%%%%%%%%%%%%%%
\section{Endoscopic transfer factor}
%%%%%%%%%%%%%%%%%%%%%%%%%%%%%%%%%%%%%%%%%%%%%%%

Many important questions in harmonic analysis on  Lie groups boil
down to the study of distributions on groups that are invariant under conjugacy.
The fundamental objects of invariant harmonic analysis are orbital integrals
as the geometric objects and characters of representations as the spectral objects.
The correspondence of these two kinds of objects reflects the core idea of harmonic analysis.

The orbital integrals are parameterized by the set of regular semisimple conjugacy classes in $G$.
Recall for such a $\gamma$, the orbital integral is defined as
$$\caO_\gamma(f)=\int_{G/G_\gamma}f(x^{-1}\gamma x) dx, \ \ f\in C^\infty_c(G),$$
and the stable  orbital integral is defined as
$$S\caO_\gamma(f)=\sum_{\gamma'\in S(\gamma)}\caO_{\gamma'}(f),$$
where $S(\gamma)$ is the stable conjugacy class.

Let $1\!\!1$ denote the trivial representation of $G$ and $\theta_{1\!\!1}$
the character of the Dirac index of the trivial representation.  That is
$$\theta_{1\!\!1} = \ch H_D^+(1\!\!1)-\ch H_D^-(1\!\!1)=\ch S^+ -\ch S^-.$$
We note that
$$\overline{\theta_{1\!\!1}} =(-1)^q  (\ch S^+ -\ch S^-)=(-1)^q \theta_{1\!\!1},$$
where $q={1\over 2} \dim G(\bbR)/K(\bbR)$.

Recall that $\theta_{\pi}$ denotes the character of the Dirac
index of $\pi$.  If $\pi$ is the discrete series representation with Dirac cohomology $E_\mu$,
then $$\theta_{\pi}= (-1)^q  \chi_\mu.$$
Labesse showed that there exists a function $f_\pi$ so that for any admissible
representations $\pi'$,
$$\tr \pi'(f_\pi)=\int_K \Theta_{\pi'}(k)\overline{\theta_{1\!\!1}\cdot \theta_{\pi}}dk.$$
Denote by $\theta_{\pi'}$ the character of its Dirac index for $\pi'$.  Then one has
$$\tr \pi'(f_\pi)=(-1)^q\int_K \theta_{\pi'}\cdot \overline{\theta_{\pi}}dk.$$

Let $\pi'$ be a discrete series representation with Dirac cohomology $E_{\mu'}$.
It follows that
$$\tr \pi'(f_\pi)=\int_K \Theta_{\pi'}(k)\overline{\theta_{1\!\!1}\cdot \theta_{\pi}}dk(\chi_{\mu'},\chi_{\mu})=\dim \Hom_K(E_{\mu'},E_{\mu}).$$
Consequently  we prove the following theorem due to Labesse.

\begin{thm}[Labesse \cite{Lab1}]
The function $f_\pi$ is a pseudo-coefficient for the discrete series $\pi$,
i.e., for any irreducible tempered representation $\pi'$,
$$\tr \pi'(f_\pi)=\begin{cases}1 \text{\ \ \ if \ }\pi\cong \pi'\\
                                                        0\text{\ \ \ \ otherwise}.\\  \end{cases} $$
\end{thm}

\begin{rmk}\label{character}
The orbital integrals of the pseudo-coefficient $f_\pi$ are easily computed for $\gamma$ regular semisimple:
$$\caO_\gamma(f_\pi)=\begin{cases}\Theta_\pi(\gamma^{-1}) \text{\ \ \ if  $\gamma$ is elliptic}\\
                                                        0\text{\ \ \ \ \ \ \ \ \ \ \  if $\gamma$ is not elliptic}.\\  \end{cases} $$
\end{rmk}

In the Langlands program a cruder form of conjugacy called stable conjugacy plays an important role.
The study of Langlands functoriality often leads to correspondence that is defined only up to stable conjugacy.
The endoscopy theory investigates the difference between ordinary and stable conjugacy
and how to understand ordinary conjugacy inside stable conjugacy.
The aim is to recover orbital integrals and characters from endoscopy groups.

The endoscopy theory for real groups is established by Shelstad in a series of papers
[Sh1-5].
Recasting Shelstad's work explicitly
in terms of the general transfer factors defined later by Langlands and
Shelstad [LS] is the first of the `Problems for Real Groups' proposed by Arthur \cite{A3}.

Recall that $G$ is a connected reductive algebraic group defined over $\bbR$.
Denote by $G^\vee$ the complex dual group and ${}^LG$ the $L$-group which is
the semidirect product of $G^\vee$ and the Weil group $W_\bbR$.  A Langlands
parameter is an $L$-homomorphism
$$\phi\colon W_\bbR\rightarrow {}^LG.$$
Two Langlands parameters are equivalent if they are conjugated by an inner automorphism of $G^\vee$.
An equivalence class of Langlands parameters is associated
to a packet of irreducible admissible representations of
$G(\bbR)$ [L2].  The $L$-packets  of Langlands parameters with bounded image
consist of tempered representations. Temperedness is respected by
$L$-packets, but not unitarity.

The discrete series $L$-packets are in bijection with the irreducible finite-dimensional
 representations of the same infinitesimal character.
 One can construct all tempered
 irreducible representations using unitary parabolic induction and by taking subrepresentations.
 Two tempered irreducible representations $\pi$ and $\pi'$ are in the same
 $L$-packet if up to equivalence, $\pi$ and $\pi'$ are subrepresentations of parabolically
 induced representations from discrete series $\sigma$ and $\sigma'$ in the same
 $L$-packets.

 A stable distribution is any element of the closure of the space spanned by all distributions
 of the form $\sum_{\pi\in \Pi}\Theta_\pi$ for $\Pi$ any tempered $L$-packet.
 Such distributions can be transferred to inner forms of $G$ via the matching
 of the stable orbital integrals, while unstable distributions cannot be.

In the setting of endoscopy embedding
$$\xi: {}^LH\rightarrow {}^LG,$$
one has a map from Langlands parameters for $H$ to that for $G$.
The Langlands functoriality principle asserts that there should be a
map from the Grothendieck group of virtual representations of $H(\bbR)$ to
that of $G(\bbR)$, compatible with $L$-packets.

We follow Labesse \S6.7 \cite{Lab2} for the description of the endoscopic transfer.
Let  $T$ be an elliptic torus of $G$ and
 $\kappa$ an endoscopic character. Let $H$ be the endoscopic
 group defined by $(T,\kappa)$.  Let $B_G$ be a Borel subgroup of $G$ containing $T$.  Set
 $$\Delta_B(\gamma)=\Pi_{\alpha>0}(1-\gamma^{-\alpha}),$$
 where the product is over the positive roots defined by $B$.  There is only one choice
 of a Borel subgroup $B_H$ in $H$, containing $T_H$ and compatible with the isomorphism
 $j\colon T_H\cong T$.

Assume $\eta\colon {}^LH\rightarrow {}^LG$ is an admissible embedding (see \S6.6 \cite{Lab2}).
Then for any pseudo-coefficent $f$ of a discrete series of $G$, there is a linear combination
$f^H$ of pseudo-coefficents of discrete series of $H$ such that for $\gamma=\gamma_G=j(\gamma_H)$ regular
in $T(\bbR)$ (see Prop. 6.7.1 \cite{Lab2}), one has
\begin{equation}\label{endoscopy}
\caS\caO_{\gamma_H}(f^H)=\Delta(\gamma_H,\gamma)\caO^\kappa_{\gamma}(f),
\end{equation}
 where the transfer factor
\begin{equation*}\label{transfer}
\Delta(\gamma_H,\gamma)=(-1)^{q(G)-q(H)}\chi_{G,H}(\gamma)\Delta_B(\gamma^{-1})\Delta_{B_H}(\gamma_H^{-1})^{-1}.
\end{equation*}

 The transfer $f\mapsto f^H$ of the pseudo-coefficents of discrete series can be extended to all of
functions in $C^\infty_c(G(\bbR))$ with extension of the correspondence $\gamma\mapsto \gamma_H$ (see
Theorem 6.7.2 \cite{Lab2}) so that the above identity (\ref{endoscopy}) holds for all $f$.

The geometric transfer $f\mapsto f^H$ is dual of a transfer for representations.
Given any admissible irreducible representation $\sigma$ of $H(\bbR)$, it corresponds to
an element $\sigma_G$ in the Grothendieck group of virtual representations of $G(\bbR)$ as follows.
Let $\phi$ be the Langlands parameter for $\sigma$.
Let $\Sigma$ be the $L$-packet of the admissible irreducible representations
of $H(\bbR)$ corresponding to a Langlands parameter $\phi$ and $\Pi$ the L-packet of representations of
$G(\bbR)$ corresponding to $\eta \circ \phi$ (that can be an empty set if this parameter
is not relevant for $G$).

\begin{thm} [Theorem 4.1.1 \cite{S}, Theorem 6.7.3 \cite{Lab2}] There is a function
$$\epsilon\colon \Pi \rightarrow \pm 1$$
such that, if we consider $\sigma_G$ in the Grothendieck group defined by
$$\sigma_G=\sum_{\pi\in\Pi}\epsilon(\pi)\pi$$
 then the transfer $\sigma\mapsto \sigma_G$ satisfies
 $$\tr \sigma_G(f)=\tr \sigma(f^H).$$
 \end{thm}

In the following we suppose that $G(\bbR)$ has a compact maximal torus $T(\bbR)$,
and $\rho-\rho_H$ the difference of half sum of positive
roots for $G$ and $H$ respectively, defines a character of $T(\bbR)$.
In  \S 7.2 of [Lab2] Labesse shows that  the canonical
 transfer factor:
 $$\Delta(\gamma^{-1})=(-1)^{q(G)-q(H)}
 \frac{\sum_{w\in W(\frg)}\epsilon (w)\gamma^{w\rho}}{\sum_{w\in W(\frh)}\epsilon (w)\gamma^{w\rho_H}}$$
 is well-defined function.
 Then the transfer factor can be expressed more explicitly if $H$ is a subgroup of $G$.
Suppose that $\frg=\frh\oplus \frs$ is the orthogonal decomposition with
respect to a non-degenerate invariant bilinear form so that the form is non-degenerate on $\frs$.
 We write $S(\frg/\frh)$ for the spin-module of the Clifford algebra $C(\frs)$.  Then
$$\Delta(\gamma^{-1})=\ch S^+(\frg/\frh) -\ch S^-(\frg/\frh).$$
In other words,  $\Delta(\gamma^{-1})$ is equal to the character of the Dirac index of
the trivial representation with respect to the Dirac operator $D(\frg,\frh)$.
 If $\Theta_\pi$ is the character of a finite-dimensional representation $\pi$,
 then
 $$\Delta(\gamma^{-1})\Theta_\pi$$ is the character of the Dirac index of $\pi$.
 This character can be calculated easily from the Kostant formula in Theorem 4.7.
 We denote by $F_\lambda$ the irreducible finite-dimensional representation of $G(\bbR)$
 with highest weight $\lambda$ and by $E_\mu$ irreducible  finite-dimensional representation
 of $H(\bbR)$ with highest weight $\mu$. Then
 $$\Delta(\gamma^{-1})\Theta_{F_\lambda}= \sum_{w\in W^1}\Theta_{E_w(\lambda+\rho)-\rho_\frh}.$$
 Here $W^1$ is a subset of elements in $W$ corresponding to $W_\frh \backslash W$ as before.

It is straightforward to use the transfer factor to calculate lifting of discrete series characters.
This lifting is closely related to the geometric transfer of the pseudo-coefficents of discrete series.
The Harish-Chandra module of a discrete series representation is isomorphic to $A_\frb(\lambda)$
for some $\theta$-stable Borel subalgebra and corresponding Harish-Chandra parameter is $\lambda+\rho$.
It follows from  \cite[Theorem 7.5]{DH} that a  tempered  representations with nonzero Dirac cohomology is
$\pi_{\lambda+\rho}=A_\frb(\lambda)$ and it has Dirac cohomology equal to an irreducble $K$-module $E_{\lambda+\rho_n}$.
In the equal rank case, it is simply a limit of discrete series.
The calculation for discrete series extends to tempered representations with nonzero Dirac cohomology.

\begin{prop} Let $\pi_\lambda$ be a discrete series of $G$ with Harish-Chandra parameter $\lambda$.
Then we have
 $$\Delta(\gamma^{-1})\Theta{\pi_\lambda}=\sum_{w\in W_{K}^1} \text{sign}(w)\Theta_{\tau_{w\lambda}}.$$
The above formula extends to limits of discrete series.
 \end{prop}
 
 \proof
 In view of Remark \ref{character}, the right hand side of (\ref{endoscopy})
 is the Dirac index of a combination of discrete series of $G(\bbR)$
 and the left hand side is a linear combination of discrete series of $H(\bbR)$.
 It follows from the Harish-Chandra formula for the character of discrete series
 and supertempered distributions \cite{HC3} that
 the Dirac index of a discrete series $\pi_\lambda$ with Harish-Chandra parameter $\lambda$ is
 $$\Delta(\gamma^{-1})\Theta{\pi_\lambda}=\sum_{w\in W_{K}^1} \text{sign}(w)\Theta_{\tau_{w\lambda}}.$$
Here ${\tau_{w\lambda}}$ denotes the discrete series for $H(\bbR)$ with Harish-Chandra parameter $w\lambda$, and
$W_{K}^1$ is a subset of elements in $W_{K}$ corresponding to $W_{H\cap K} \backslash W_{K}$.
This calculation is compatible with Labesse's calculation of the transfer of the pseudo-coefficients
of discrete series in \S 7.2 \cite{Lab2}.

It remians to show that this formula extends to  tempered elliptic representations nonzero Dirac cohomology.  By 
Theorem 7.5 of \cite{DH} that a  tempered elliptic representations nonzero Dirac cohomology is a limit of discrete sereis
$\pi_{\lambda+\rho}=A_\frb(\lambda)$ and it has Dirac cohomology equal to an irreducble $K$-module $E_{\lambda+\rho_n}$.
As the parmeter for a limit of discrete series,  $\lambda+\rho$ is regular with respect to compact roots.
Thus, the same calculation for discrete series applies here. 
\qed

\medskip

As a final remark, we note that  we may use the Arthur packets to deal with  non-tempered case.
 The Arthur packets are parameterized by mappings
$$\psi\colon  W_\bbR \times SL(2,\bbC) \rightarrow {}^LG$$
for which the projection onto the dual group $G^\vee$ of $\psi(W_\bbR)$ is relatively compact.
Adams and Johnson \cite{AJ} have constructed some $A$-packets consisting of unitary $A_\frq(\lambda)$-modules.
As most of unitary $A_\frq(\lambda)$-modules can be classified
by their Dirac cohomology \cite{HPV}, the determination of Dirac cohomology of
$A_\frq(\lambda)$-modules may have some bearing on answering Arthur's questions (See \cite[\S9]{A2})
on Arthur packet $\Pi_\psi$.


\begin{thebibliography}{22}



\bibitem{A} J. Adams,
\emph{Lifting of characters on orthogonal and metaplectic groups},
Duke Math. J. \textbf{92} (1998), 129--178.

\bibitem{AJ} J. Adams and J. Johnson,
\emph{Endoscopic groups and packets of non-tempered representations},
Compos. Math. \textbf{64} (1987), 271--309.

\bibitem {A1} J. Arthur, \emph{On elliptic tempered characters},
Acta Math. \textbf{171} (1993), 73--138.

\bibitem{A2} J. Arthur, \emph{On the Fourier transforms of weighted orbital integrals},
J. Reine Angew. Math. \textbf{452} (1994), 163--217.

\bibitem{A3} J. Arthur, \emph{Problems for real groups},
in book {Representation Theory of Real Reductive Lie Groups}, Contemporary Mathematics
\textbf{472} (2008), 39--62.

\bibitem{AS} M.~Atiyah, W.~Schmid, \emph{A geometric
construction of the discrete series for semisimple Lie groups}
Invent. Math. \textbf{42} (1977), 1-62.

\bibitem{AS'} M.~Atiyah, W.~Schmid,
\emph{Erratum: A geometric construction of the discrete series for
semisimple Lie groups}, Invent. Math. \textbf{54} (1979), 189--192.

\bibitem{DH} C.-P. Dong, J.-S. Huang, \emph{Dirac cohomology of cohomologically induced modules},
Amer. J. Math.  {\bf 137} (2015), 37-60.

\bibitem{HC1} Harish-Chandra, \emph{The characters of semisimple Lie groups}, Trans. Amer.
Math. Soc. \textbf{83} (1956), 98-163.

\bibitem{HC2} Harish-Chandra, \emph{Harmonic analysis on real reductive groups I.
The theory of the constant term}, J. Funct. Anal. \textbf{19} (1975),
104-204.

\bibitem{HC} Harish-Chandra,
\emph{Discrete series for semisimple Lie groups, I and II}, Acta
Math. \textbf{113} (1965), 242-318 and \textbf{116} (1966), 1-111.

\bibitem{HC3} Harish-Chandra, \emph{Supertempered Distributions on Real Reductive Groups},
Studies in Applied Mathematics,  Advances in Mathematics Supplementary Studies, \textbf{8} (1983),
139-153. (It is contained in Harish-Chandra's collected Collected Papers, Volumed IV, 447-461.)

\bibitem{Hd} T.~Howard, \emph{Lifting of characters on p-adic orthogonal and metaplectic groups}, 
Compos. Math.  \textbf{146} (2010), 795-810.

\bibitem{Ho} R.~Howe, \emph{Remarks on classical invariant
theory}, Trans. Amer. Math. Soc. \textbf{313} (1989), no. 2, 539--570.


\bibitem{H} J.-S.~Huang, \emph{Dirac cohomology, elliptic representations and endoscopy},
Representation Theory of Reductive Groups, in Honor of 60th Birthday of David Vogan,
Birkhauser, Progress in Mathematics, {\bf 321} (2015), 241-276.

\bibitem{H2} J.-S.~Huang, \emph{Symplectic Dirac cohomology and lifting of characters to metaplectic groups},
arXiv:2006.00157.



\bibitem{HKP} J.-S.~Huang, Y.-F.~Kang and P.~Pand\v zi\'c, \emph{Dirac
cohomology of some Harish-Chandra modules}, Transform. Groups \textbf{14} (2009), no. 1, 163-173.

%\bibitem{HPR} J.-S Huang, P. Pand\v{z}i\'{c}, D. Renard, \emph{Dirac
%operators and Lie algebra cohomology}, Represent. Theory \textbf{10} (2006), 299--313.

\bibitem{HP1} J.-S Huang, P. Pand\v{z}i\'{c}, \emph{Dirac
cohomology, unitary representations and a proof of a conjecture of
Vogan}, J. Amer. Math. Soc. \textbf{15} (2002), 185--202.

\bibitem{HP2} J.-S Huang, P. Pand\v{z}i\'{c}, \emph{Dirac operator in
Representation Theory}, Math. Theory Appl., Birkh\"{a}user, 2006.

\bibitem{HP3} J.-S Huang, P. Pand\v{z}i\'{c}, \emph{Dirac
cohomology for Lie superalgebras}, Transform. Groups \textbf{10} (2005), 201--209.

\bibitem{HPV} J.-S.~Huang, P.~Pand\v zi\'c and D.-A.~Vogan,
\emph{On classifying unitary representations by their Dirac cohomology},
Science China Mathematics \textbf{60} (2017), 1937-1962.


\bibitem{HPZ} J.-S.~Huang, P.~Pand\v zi\'c and F.-H.~Zhu,
\emph{Dirac cohomology, K-characters and branching laws},
Amer.  J. Math. \textbf{135} (2013), 1253-1269.






\bibitem{Ko1} B. Kostant, \emph{A cubic Dirac operator and the
emergence of Euler number multiplets of representations for equal
rank subgroups}, Duke Math. J. \textbf{100} (1999), 447--501.



\bibitem{Ko2} B. Kostant, 
\emph{Dirac cohomology for the cubic Dirac
operator, Studies in Memory of Issai Schur}, Progress in Math. Vol.
\textbf{210} (2003), 69--93.


\bibitem{Lab1} J.-P. Labesse, 
\emph{Pseudo-coefficients tr\`es cuspidaux et K-th\'eorie},
Math. Ann \textbf{291} (1991), 607-616.

\bibitem{Lab2} J.-P. Labesse, \emph{Introduction to endoscopy}, in
Representation Theory of Real Reductive Lie Groups, Contemporary Mathematics \textbf{472} (2008),
197-213.


\bibitem{L2} R.~P. Langlands, 
\emph{On the classification of irreducible representations of real algebraic groups},  Math. Surveys and Monographs \textbf{31} (1989), AMS, 101-170.

\bibitem]{LS} R.~P. Langlands, D. Shelstad,
 \emph{On the definition of transfer factors},  Math. Ann. \textbf{278} (1987), 219--271.


\bibitem{Li} W.-W. Li, \emph{Spectral transfer for metaplectic groups. I. Local character relations},  J. Inst. Math. Jussieu \textbf{18} (2019), 25-123.


\bibitem{P} R.~Parthasarathy, \emph{Dirac operator and the discrete
series}, Ann. of Math. \textbf{96} (1972), 1--30.


\bibitem{R} D. Renard,
\emph{Endoscopy for $Mp(2n,\bbR)$},  Amer. J. Math. \textbf{121} (1999), 1215-1243.


\bibitem{SR} S.~A. Salamanca-Riba,
\emph{On the unitary dual of real reductive Lie groups
and the $A_\frq(\lambda)$ modules: the strongly regular case},
 Duke Math. J. \textbf{96} (1998), 521-546.

\bibitem{Sh1} D. Shelstad, \emph{Notes on L-indistinguishability (based on a lecture by R.P. Langlands)}, in
Automorphic Forms, Representations and L-functions, Proc. Sympos. Pure Math.
\textbf{33}, Part 2, Amer. Math. Soc., 1979, 193--204.

\bibitem{Sh2} D. Shelstad, \emph{Characters and inner forms of a quasi-split group over R},  Compositio
Math. \textbf{39} (1979), 11--45.

\bibitem{Sh3} D. Shelstad, \emph{Orbital integrals and a family of groups attached to a real reductive group},
Ann. Sci. Ecole Norm. Sup. \textbf{12} (1979), 1--31.

\bibitem{Sh4}  D. Shelstad, \emph{Embeddings of L-groups}, Canad. J. Math. \textbf{33} (1981), 513--558.

\bibitem{S} D. Shelstad, \emph{L-indistinguishability for real groups}, Math. Ann. \textbf{259} (1982), 385--430.


\bibitem{V2} D.~A.~Vogan, Jr., \emph{Irreducible characters of semisimple
Lie groups II.  The Kazhdan-Lusztig conjectures}, Duke Math. J.
\textbf{46} (1979), 805--859.

\bibitem{V1} D.~A.  Vogan Jr., \emph{Representations of Real Reductive
Lie Groups}, Progress in Math. Vol. \textbf{15}, Birkh\"{a}user, 1981.

\bibitem{V} D.~A. Vogan Jr., \emph{Dirac operators and unitary
representations}, 3 talks at MIT Lie groups seminar, Fall 1997.

\bibitem{VZ} D.~A. Vogan, Jr., G. J. Zuckerman,
\emph{Unitary representations with non-zero cohomology}, Compos. Math. \textbf{53} (1984), 51--90.


\bibitem{W} J.~A.~Wolf, \emph{Partially harmonic spinors and representations of reductive Lie groups}, J. Funct. Anal.
\textbf{15} (1974), 117--154.








\end{thebibliography}
\end{document}